\documentclass[11pt,twoside]{article}
\def\wdeb{0}
\usepackage{amsfonts}
\usepackage{bezier}
\usepackage{euscript}
\usepackage{amsfonts,amssymb,amsmath,amsbsy,amsthm,amscd}
\usepackage[cp1251]{inputenc}
\usepackage[english,russian]{babel}   
\usepackage{allerree}
\begin{document}

\cleanbegin

\def\udk{514.17+514.174}
\ltitle{Оценка расстояния между двумя телами внутри $n$-мерного единичного куба и шара}
{Ф.\,А.~Ивлев\footnote[1]{%
{\it Федор Алексеевич Ивлев}~--- студ. каф. высшей алгебры мех.-мат. ф-та МГУ, e-mail: ivlevfyodor@gmail.com.}
}

\iabstract{
Рассматривается задача об оценке расстояний между двумя телами объема~$\varepsilon$, расположенными внутри $n$-мерного тела~$B$ единичного объема, при $n \to \infty$.
В ряде случаев такие расстояния ограничены функцией от $\varepsilon$, не зависящей от~$n$.
Мы рассматриваем случаи, когда $B$~--- шар или куб.
}
{
Минимальная поверхность, многомерная выпуклая геометрия, центральные предельные теоремы.
}
{
The problem of bounding of the distance between the two bodies of volume~$\varepsilon$ located inside the $n$-dimensional body~$B$ of unit volume where $n \to \infty$ is considered.
In some cases such distances are bounded by function depends on~$\varepsilon$ but not depends on~$n$.
We consider cases when $B$ is a sphere or a cube.}
{
Minimal surface, multidimensional convex geometry, central limit theorems.
}


\textbf{1. Введение}

Пусть $B_n$~--- тело единичного объема в $n$-мерном пространстве.
При $n \to \infty$ диаметр~$B_n$ тоже стремится к бесконечности,  так что внутри тела~$B$ можно найти далекие точки.
Это обстоятельство связано с  трудностями при переносе конечномерных результатов на бесконечномерные, в частности, в функциональном анализе.
Тем не менее в ряде случае есть основания предполагать, что если взять два множества объема~$\varepsilon$, то расстояние между ними окажется ограниченной функций от $\varepsilon$, вне зависимости от~$n$. Это обстоятельство может оказаться полезным, в том числе и для переноса результатов на бесконечномерный случай. Интересно, что результаты затрагивают теорию минимальных поверхностей [2], [3].

Под {\it расстоянием} $\mathrm{dist}(A, B)$ между множествами $A$ и $B$ мы понимаем величину:

$$
D = \mathrm{dist}(A, B) = \inf_{X \in A, Y \in B} \mathrm{dist}(X, Y),
$$

где $\mathrm{dist}(X, Y)$ есть расстояние между точками $X$ и $Y$.

Следующие две гипотезы были предложены Н.~А.~Бобылевым и А.\,Я.\,Канелем: 

\begin{theorem}{Гипотеза 1. (случай куба)}
Пусть $\varepsilon$~--- данное число в интервале $(0,1)$, $K_n$~--- $n$-мерный куб единичного объема.
Внутри $K_n$ выбраны два множества $A$ и $B$, каждое объема $\varepsilon$.
Тогда расстояние между $A$ и $B$ не больше, чем некоторая константа $D = D(\varepsilon)$, и не зависит от размерности пространства~$n$.
\end{theorem}

\begin{theorem}{Гипотеза 2. (случай шара)}
Пусть $\varepsilon$~--- данное число в интервале $(0,1)$, $K_n$~--- $n$-мерный шар единичного объема.
Внутри $K_n$ выбраны два множества $A$ и $B$, каждое объема $\varepsilon$.
Тогда расстояние между $A$ и $B$ не больше, чем некоторая константа $D = D(\varepsilon)$, и не зависит от размерности пространства~$n$.
\end{theorem}

Очевидно, что, если $\varepsilon \geq 1/2$, то обе гипотезы верны.
Поэтому в дальнейшем мы будем предполагать, что  $\varepsilon < 1/2$.

В данной работе гипотеза 1 доказана в случае, когда множества~$A$ и $B$ являются пересечениями полупространств с кубом, причем границы этих полупространств перпендикулярны главной диагонали куба.
Также доказана гипотеза 2 в случае, когда множества~$A$ и $B$ выпуклы.
В конце приведены общая идея о том, почему поставленная гипотеза скорее всего верна и для произвольных множеств, план доказательства гипотезы, постановки близких задач.


\textbf{2. Выпуклые множества}

Зачастую для удобства мы будем опускать индекс размерности у рассматриваемых тел и, например, обозначать $n$-мерный куб просто $K$.

\textbf{2.1. Общий случай}

Сделаем несколько общих замечаний, применимых к поставленной задаче для любой фигуры.

Вместо множеств $A$ и $B$ можно взять из замыкания $\overline{A}$ и $\overline{B}$ и расстояние от этого между ними не изменится.
Поэтому в дальнейшем мы будем считать, что множества $A$ и $B$ замкнуты, и, следовательно, компактны.

Напомним классический факт:

\begin{theorem}{Лемма 1.}
	Пусть $A$ и $B$ выпуклые компактные подмножества $\mathbb{R}^n$.
	Тогда существуют параллельные гиперплоскости $\Pi_A$ и $\Pi_B$, разделяющие множества $A$ и $B$, расстояние между которыми равно расстоянию между множествами $A$ и $B$.
\end{theorem}

Обозначим рассматриваемую фигуру единичного объема через $F$.
Рассмотрим часть~$F$, отсекаемую гиперплоскостью $\Pi_A$, и содержащую множество~$A$.
Обозначим ее через~$A'$.
Аналогично определим~$B'$.
Тогда множества~$A'$ и $B'$ объема хотя бы~$\varepsilon$, потому что содержат одно из множеств~$A$ и $B$, а расстояние между ними такое же, как расстояние между~$A$ и $B$.
Если заменить в гипотезе условие, что объемы множества~$A$ и $B$ равны~$\varepsilon$, на то, что объемы этих множеств не меньше $\varepsilon$, мы получим гипотезу очевидно равносильную исходной.
Поэтому при работе с выпуклыми множествами мы сразу будем оперировать с множествами вида $A'$ и $B'$ и называть просто $A$ и $B$ соответственно.

\textbf{2.2. Случай шара}

\begin{theorem}{Теорема 1.}
	Пусть задано число
	$\varepsilon \in (0, 1/2)$,
	$S_n$~--- $n$-мерный шар единичного объема.
	Выпуклые подмножества $A$ и $B$ этого шара имеют объем~$\varepsilon$ каждое.
	Тогда расстояние между $A$ и $B$ не превосходит некоторой константы~$D = D(\varepsilon)$, не зависящей от~$n$.
\end{theorem}

\textbf{Доказательство}
	В силу леммы 1 имеем гиперплоскости~$\Pi_A$ и $\Pi_B$, и множества~$A$ и $B$ суть пересечения одного из полупространств, на которые делят пространство гиперплоскости~$\Pi_A$ и $\Pi_B$ соответственно, с шаром~$S$.
	Будем оценивать расстояние между $A$ и $B$.
	
	Обозначим центр шара через~$O$, а плоскость, проходящую через~$O$ параллельную~$\Pi_A$ через~$\Pi$.
	Искомое расстояние равно расстоянию между параллельными, как следует из леммы 1, плоскостями~$\Pi_A$ и $\Pi_B$.
	Значит, оно равно удвоенному расстоянию от точки~$O$ до множества~$A$ или, что то же самое, расстоянию между плоскостями~$\Pi_A$ и $\Pi$.
	Обозначим это расстояние через~$d$.
	
	\begin{theorem}{Лемма 2.}
		Объем части шара $S$ единичного объема, находящейся между гиперплоскостью~$\Pi$, проходящей через его центр и параллельной гиперплоскостью~$\Pi'$, отстоящей от первой на расстояние~$d$, при стремлении размерности~$n$ к бесконечности стремится к
		$$
		\sqrt{e} \int_0^d{e^{-\pi e x^2}\,dx}
		$$
	\end{theorem}
	
	\textbf{Доказательство.}
		Очевидно этот объем равен
		$$
		V = \int_0^d{S(x)dx},
		$$
		где $S(x)$~--- объем $(n-1)$-мерного шара, являющегося сечением исходного шара~$S$ гиперплоскостью~$\Pi_x$, параллельной~$\Pi$ и находящейся от нее на расстоянии~$x$.
		Ее радиус по теореме Пифагора равен $r = \sqrt{R_n^2 - x^2}$, где $R_n$~--- радиус $n$-мерного шара единичного объема.
		Объем этого шара равен $C_{n-1} r^{n-1}$, где
		$C_{n-1} = \dfrac{\pi^{(n-1)/2}}{\Gamma \left(\frac{n - 1}{2} + 1\right)}$.
		Из уравнения $C_n R_n^n = 1$ находим формулу для $R_n$
		$$
		R_n = \frac{\Gamma\left(\frac{n}{2} + 1\right)^{1/n}}{\pi^{1/2}}.
		$$
		Подставляя это выражение в формулу объема имеем:
		\begin{gather*}
		V = \int_0^d{S(x)dx} =
		\int_0^d{C_{n-1}\left(\sqrt{R_n^2 - x^2}\right)^{n - 1}\,dx} =\\
		= \int_0^d{C_{n-1}R_{n-1}^{n-1} \cdot \left(\frac{R_n}{R_{n-1}}\right)^{n-1} \cdot \left(\sqrt{1 - \frac{x^2}{R_n^2}}\right)^{n-1}\,dx} =
		\left(\frac{R_n}{R_{n-1}}\right)^{n-1} \int_0^d{\left(\sqrt{1 - \frac{x^2}{R_n^2}}\right)^{n-1}\,dx}.
		\end{gather*}
		
		Найдем к чему стремится первый множитель.
		Имеем
		$$
		\left(\dfrac{R_n}{R_{n-1}}\right)^{n-1} =
		\left(
		\dfrac{\Gamma\left(\frac{n}{2} + 1\right)^{1/n}
			}{
			\Gamma\left(\frac{n-1}{2} + 1\right)^{1/(n-1)}
			}
		\right)^{n-1} =
		\dfrac{\Gamma \left(\frac{n}{2} + 1\right)^{\tfrac{n - 1}{n}}}
			{\Gamma \left(\frac{n - 1}{2} + 1\right)}
		$$
		Воспользуемся асимптотической формулой роста Гамма-функции в следующем виде
		$$
			\Gamma(z + 1) =
			e^{-z} z^{z + 1/2} (2\pi)^{1/2} (1 + O\!\left(\tfrac{1}{z}\right)), \quad \text{при } z \to \infty.
		$$
		Получаем в нашем случае
		\begin{gather*}
			\dfrac{\Gamma \left(\frac{n}{2} + 1\right)^{\tfrac{n - 1}{n}}}
				{\Gamma \left(\frac{n - 1}{2} + 1\right)} =
			\dfrac{
				\left(e^{-n/2} \left(\frac{n}{2}\right)^{(n + 1) / 2} (2 \pi)^{1/2} (1 + O(\frac{1}{n}))\right)^{\tfrac{n - 1}{n}}
				}{
				e^{-(n - 1) / 2} \left(\frac{n - 1}{2}\right)^{n / 2} (2 \pi)^{1 / 2} (1 + O(\frac{1}{n}))
				}
			=\\
			= \dfrac{
				e^{-(n - 1) / 2} \left(\frac{n}{2}\right)^{n / 2} \left(\frac{n}{2}\right)^{\frac{-1}{2n}} (2 \pi)^{1 / 2} (2 \pi)^{\frac{-1}{2n}} (1 + O(\frac{1}{n}))^{\frac{n - 1}{n}}
				}{
				e^{-(n - 1) / 2} \left(\frac{n - 1}{2}\right)^{n / 2} (2 \pi)^{1 / 2} (1 + O(\frac{1}{n}))
				}
			=
			\left(\dfrac{n}{n - 1}\right)^{\tfrac{n}{2}} (2 \pi)^{\tfrac{-1}{2n}} \dfrac{(1 + O(\frac{1}{n}))^\frac{n-1}{n}}{1 + O(\frac{1}{n})}
			=\\
			= \left(1 - \frac{1}{n - 1}\right)^{\tfrac{n - 1}{2}} (1 + (1)) = \exp(\tfrac{1}{2}) (1 + (1)), \quad \text{при } n \to \infty.
		\end{gather*}
		
		Для доказательства леммы осталось показать, что
		$$
			\lim_{n \to \infty}{\int_0^d{\left(\sqrt{1 - \dfrac{x^2}{R_n^2}}\right)^{n-1}dx}} =
			\dfrac{1}{\sqrt{\pi e}} \int_0^d{e^{-x^2}dx}.
		$$
		Обозначим подынтегральное выражение через $F(x, n)$.
		Заметим, что оно положительно начиная с некоторого $n$ и не превосходит единицы.
		Поэтому

		Покажем, что последовательность $F(x, n)$ равномерно по $n$ на отрезке $[\delta, d]$ сходится к $\exp(-x^2 \pi e)$.
		Для этого представим $F(x, n)$ в следующем виде:
		$$
			F(x, n) =
			\left(1 - \dfrac{1}{R_n^2/x^2}\right)^{\displaystyle R_n^2/x^2 \cdot \frac{n-1}{2(R_n^2/x^2)}} =
			\left(\left(1 - \dfrac{1}{R_n^2/x^2}\right)^{\displaystyle R_n^2/x^2}\right)^{\displaystyle \frac{n-1}{2(R_n^2/x^2)}}.
		$$
		Поскольку $\lim_{y \to \infty}{(1 - 1/y)^y} = e^{-1}$, для любого $\alpha > 0$ существует такое $Y > 1$, что для любого $y > Y$ верно, что
		$$
			(1 - 1/y)^y = e^{-1 + \beta}, \quad \text{где } |\beta| < \alpha.
		$$
		
		Так как $R_n \to \infty$ при $n \to \infty$, то начиная с некоторого $N_1$ для всех $n > N_1$ верно, что $R_n > dY$.
		Следовательно, для всех $x \in [\delta, d]$ верно, что $R_n^2/x^2 > Y$, а значит,
		$$
			\left(1 - \dfrac{1}{R_n^2/x^2}\right)^{R_n^2/x^2} =
			e^{-1 + \beta(n)}, \quad \text{где } |\beta(n)| < \alpha \ \text{для всех} n > N_1.
		$$
		Тогда для $n > N_1$ имеем
		$$
			F(x, n) =
			\left(e^{-1 + \beta(n)}\right)^\frac{x^2(n-1)}{2 R_n^2} =
			\left(e^{(-1 + \beta(n))((n-1)/(2R_n^2))}\right)^{x^2}.
		$$
		
		Заметим, что
		\begin{gather*}
			\dfrac{n - 1}{2 R_n^2} =
			\dfrac{n - 1}{2\left(\frac{\Gamma\left(\frac{n}{2} + 1\right)^{1/n}}{\pi^{1/2}}\right)^2} =
			\dfrac{n - 1}{2\left(\left(e^{-n/2} \left(\frac{n}{2}\right)^{(n + 1) / 2} (2 \pi)^{1/2} (1 + O(\frac{1}{n}))\right)^{1/n}/\pi^{1/2}\right)^2} = \\
			= \dfrac{(n - 1)\pi}{2e^{-1} \left(\frac{n}{2}\right)^{(n+1)/n} (2\pi)^{1/n} (1 + o(1))} =
			\dfrac{\pi e (n - 1)}{2\dfrac{n}{2} (1 + o(1))} =
			\pi e (1 + o(1)).
		\end{gather*}
		Следовательно, всегда можно подобрать такое $N$, что при $n > N$ функция $F(x, n)$ равномерно по $n$ для всех $x$ из отрезка $[\delta, d]$ приближает функцию $e^{-\pi e x^2}$ с наперед заданной точностью.
		
		Получаем, что
		$$
			\lim_{n \to \infty}{\int_0^d{\left(\sqrt{1 - \dfrac{x^2}{R_n^2}}\right)^{n-1}dx}} =
			\int_0^d{e^{-\pi e x^2}\,dx},
		$$
		что и требовалось.
	
	Лемма 2 доказана.
	
	Для завершения доказательства теоремы осталось заметить, что, подбирая соответствующее $d$, мы сможем получить любой наперед заданный объем $1/2 - \varepsilon \in (0, 1/2)$.
	Действительно, если, отступая от центра шара на расстояние $d$, мы набираем объем хотя бы $1/2 - \varepsilon$ начиная с некоторой размерности, то расстояние от $O$ до $A$ не превышает этого расстояния $d$.
	Следовательно, расстояние между $A$ и $B$ не будет превышать $2d$ начиная с некоторой размерности, что и требуется доказать.
	
	Покажем, что мы действительно можем, отступая на расстояние $d$, получить любой наперед заданный объем из интервала $(0, 1/2)$.
	Для этого осталось заметить, что
	$$
		\sqrt{e} \int_0^\infty{e^{-\pi e x^2}\,dx} =
		\sqrt{e} \frac{1}{\sqrt{\pi e}} \int_0^\infty{e^{-u^2}\,du} =
		\frac{1}{\sqrt{\pi}} \cdot \frac{\sqrt{\pi}}{2} =
		1/2.
	$$
	Следовательно, мы видимо, что надо расстояние между множествами будет стремиться к такому числу $D = D(\varepsilon)$, что $\sqrt{e} \in_0^D{e^{\pi e x^2}\,dx = 1/2 - \varepsilon}$.
	
	Теорема 1 доказана.

\textbf{2.3. Случай куба}

Мы опять же ограничимся рассмотрением выпуклых множеств $A$ и $B$.
Поэтому можно считать, что эти множества суть пересечения некоторых полупространств с кубом, причем границы соответствующих полупространств, гиперплоскости $\Pi_A$ и $\Pi_B$ соответственно, параллельны.
Так как нас интересует возможный максимум расстояний, то мы считаем, что объемы $A$ и $B$ в точности равны $\varepsilon$.
Значит, по направлению нормали к разделяющим гиперплоскостям можно однозначно восстановить сами гиперплоскости, множества и расстояние между ними.
Поэтому можно считать, что искомое расстояние есть функция на единичной $(n-1)$ мерной сфере.
Обозначим эту функцию через $f(\vec{n})$.
Будем считать, что координаты всех вершин куба равны либо 0 либо 1.
Тогда очевидно, что все квадранты для функции $f$ равноправны, поэтому мы будем рассматривать только квадрант, где все координаты положительны.
Так как функция определена на компакте, то она достигает своего максимума на нем.

Если какая-то координата $n_i$ нормали $\vec{n}$ равна нулю, то соответствующая ей гиперплоскость перпендикулярна любой гиперплоскости вида $x_i = \mathrm{const}$.
Но тогда в сечении гиперплоскостью такого вида получится $(n-1)$-мерный куб с множествами $A$ и $B$ $(n-1)$-мерного объема $\varepsilon$ и таким же расстоянием между ними.
Поэтому все нулевые координаты можно исключить из рассмотрения, перейдя к меньшей размерности.

Наше доказательство будет опираться на следующий недоказанный факт.

\begin{theorem}{Гипотеза 3.}
	В точке достижения глобального максимума функции $f(\vec{n})$ все ненулевые координаты нормали $\vec{n}$ равны между собой.
\end{theorem}

По сути этот факт нам говорит о том, что лучший ответ для любого $\varepsilon$ достигается в некоторой размерности, как гиперплоскость, перпендикулярная главной диагонали куба.
То есть оба множества максимально <<вжимаются>> в противоположные углы куба.
в предположении справедливости гипотезы докажем следующую теорему.

\begin{theorem}{Теорема 2.}
	Пусть задано число $\varepsilon \in (0, 1/2)$, $K_n$~--- $n$-мерный куб единичного объема.
	Выпуклые подмножества $A$ и $B$ куба $K_n$ имеют объем $\varepsilon$ каждое.
	Тогда расстояние между $A$ и $B$ не превосходит некоторой константы $d = d(\varepsilon)$, не зависящей от $n$.
\end{theorem}

\textbf{Доказательство.}
	Будем считать, что координаты всех вершин куба равны либо 0 либо 1.
	Для каждого $n$ рассмотрим плоскость $\Pi_A^n$ вида $x_1 + x_2 + \ldots + x_n = a(n)$ и соответствующее ей множество $A_n$ вида $x_1 + x_2 + \ldots + x_n \leq a(n)$, $x_i \geq 0$ для всех $i \in \{1, \ldots, n\}$.
	При этом $a(n)$ выбирается так, чтобы объем $A_n$ был равен $\varepsilon$.
	Тогда соответствующее множество $B_n$ будет симметрично $A_n$ относительно центра куба.
	Расстояние между множествами в этом случае будет равно расстоянию между проекциями этих множеств на главную диагональ куба, соединяющую вершины $(0, 0, \ldots, 0)$ и $(1, 1, \ldots, 1)$.
	Расстояние от начала координат вдоль этой диагонали пропорционально сумме координат точки с коэффициентом $1/\sqrt{n}$.
	Значит, расстояние между этими множествами будет равно
	${n \cdot 1 / \sqrt{n} - 2 \cdot a(n) \cdot 1 / \sqrt{n}} = {\sqrt{n} - 2a / \sqrt{n}}$.
	
	Для того, чтобы оценить это расстояние введем $n$ независимых равномерно распределенных на отрезке $[0, 1]$ случайных величин $\xi_1$, $\xi_2$, \ldots, $\xi_n$.
	Тогда набор значений этих случайных величин задает точку в нашем единичном кубе.
	Так как все они независимые и равномерно распределены на $[0, 1]$, то распределение внутри куба будет однородное.
	В этом случае вероятность попадания точки в наше множество будет равна его объему, то есть $\varepsilon$.
	Значит,
	\begin{gather*}
		\varepsilon =
		P\left(\sum\limits_{i = 1}^n{\xi_i} \leq a(n)\right) =
		P\left(\sum \xi_i - \sum E_{\xi_i} \leq a(n) - \sum E_{\xi_i}\right) =\\
		= P\left(\dfrac{\sum \xi_i - \sum E_{\xi_i}}{D_{\xi_1}\sqrt{n}} \leq \dfrac{a(n) - \frac{1}{2}n}{D_{\xi_1}\sqrt{n}}\right) =
		P\left(\dfrac{\sum \xi_i - \sum E_{\xi_i}}{D_{\xi_1}\sqrt{n}} \leq \dfrac{a(n) - \frac{1}{2}n}{\frac{1}{12}\sqrt{n}}\right) =\\
		= {\Phi\left(\dfrac{a(n) - \frac{1}{2}n}{\frac{1}{12}\sqrt{n}}\right) + \frac{c(n)}{\sqrt{n}}},
	\end{gather*}
	
	По теореме Берри--Эссеена\footnote{Подробное описание приведено в [1]}, примененной к нашим случайным величинами, имеем
	$$
		\varepsilon =
		P\left(\dfrac{\sum \xi_i - \sum E_{\xi_i}}{D_{\xi_1}\sqrt{n}} \leq \dfrac{a(n) - \frac{1}{2}n}{\frac{1}{12}\sqrt{n}}\right) =
		{\Phi\left(\dfrac{a(n) - \frac{1}{2}n}{\frac{1}{12}\sqrt{n}}\right) + \frac{c(n)}{\sqrt{n}}},
	$$
	где $c(n)$~--- некоторая ограниченная функция.
	
	Заметим, что в левой части этого равенства стоит константа, а второе слагаемое в правой части стремится к нулю при стремлении $n$ к бесконечности.
	Значит, первое слагаемое в правой части этого равенства стремится к $\varepsilon$.
	А следовательно,
	$$
		\lim_{n \to \infty}{\dfrac{a(n) - \frac{1}{2}n}{\frac{1}{12}\sqrt{n}}}= \Phi^{-1}(\varepsilon) = b.
	$$
	Имеем
	\begin{gather*}
		\dfrac{a(n) - \frac{1}{2}n}{\frac{1}{12}\sqrt{n}} = b + (1)\\
		a(n) = \frac{n}{2} + \frac{1}{12}b\sqrt{n} + (\sqrt{n}).
	\end{gather*}
	
	Осталось вспомнить, что искомое расстояние равно $\sqrt{n} - 2a / \sqrt{n}$.
	Обозначим его через~$d(n)$ и подставим полученное значение для $a(n)$:
	\begin{gather*}
		d(n) = \sqrt{n} - \dfrac{2(\frac{n}{2} + \frac{1}{12}b\sqrt{n} + (\sqrt{n})}{\sqrt{n}} =
		\sqrt{n} - \frac{n + \frac{1}{6}b\sqrt{n} + (\sqrt{n})}{\sqrt{n}} =
		\sqrt{n} - \sqrt{n} - \frac{1}{6}b + (1)\\
		d(n) = -\frac{1}{6}b + (1) = -\frac{1}{6}\cdot\Phi^{-1}(\varepsilon) + (1)
	\end{gather*}
	
	Получаем, что искомое расстояние стремится к $-1 / 6 \cdot\Phi^{-1}(\varepsilon)$ при $n \to \infty$.
	Значит, расстояние $d(n)$ ограничено некоторой константой, независящей от $n$, ч.\,т.\,д.
	
	Теорема 2 доказана.

\textbf{3. Общий случай}

Выше мы рассмотрели случаи выпуклых множеств $A$ и $B$ и достаточно сильно пользовались спецификой выпуклости.
В общем случае, когда $A$ и $B$ могут быть невыпуклыми, их уже нельзя будет разделить гиперплоскостями.
Даже, если это нам удастся, расстояние между этими гиперплоскостями может быть не равно расстоянию между множествами.
Автор считает верным предположение о том, что наибольшее расстояние между множествами $A$ и $B$ все-таки достигается тогда, когда они оба выпуклы.
Автору также пока неизвестно строгое доказательства этого предположения.

Можно считать, что для множеств $A$ и $B$ определена их граница, а у границы определена площадь.
Иначе $A$ и $B$ можно приблизить множествами с указанными свойствами и, перейдя к пределу, получить то же самое расстояние между множествами в пределе.
Определим площадь свободной поверхности множеств $A$ и $B$ как площадь той части их границы, которая не является границей объемлющей фигуры (куба, шара, \ldots).

Идея состоит в том, что даже при маленьком объеме, но как-то ограниченного снизу, при достаточно большой размерности пространства площадь поверхности границы этого множества будет достаточно большой.
Так же хочется показать, что вместе с ней будет достаточно большой площадь свободной поверхности.
Тогда, пользуясь тем, что объем $\delta$-окрестности множества, отличается от объема исходного множества примерно на произведение $\delta$ и площади свободной поверхности, мы получим, что объем $d$-окрестности множества $A$ будет иметь объем хотя бы $\int_0^dS(x, \varepsilon) dx$.
Здесь $S(x, \varepsilon)$ обозначает минимальную возможную площадь свободной поверхности среди всевозможных $x$-окрестностей множеств объема $\varepsilon$.

Если оценить $S(x, \varepsilon)$ снизу, то мы получим, что для некоторого фиксированного ${d = d(\varepsilon)}$ объем $d$-окрестности любого множества объема $\varepsilon$ будет равен хотя бы $1 / 2$.
Но тогда расстоянием между множествами $A$ и $B$ будет не более, чем $2d$.
Действительно, при взятии $d$- окрестности каждого из них мы получим множества объема равного хотя бы половине всего объема объемлющей фигуры.
Значит, они пересекаются и мы можем найти точки $A$ и $B$ на расстоянии не более, чем $d$ от этой точки пересечения.

\textbf{4. Возможные обобщения}

В данной работе были рассмотрены случаи подмножеств куба и шара единичного объема.
Можно обобщить утверждение основной гипотезы на случай произвольного тела.
Задача перестает быть интересной, если разрешить телу быть <<вытянутым>> относительно ограниченного числа координат:
если в качестве тела позволить брать параллелепипед вытянутый вдоль одной координаты и узкий относительно других\footnote{Например, со сторонами $0.5, 0.5, \ldots, 0.5, 2^{n-1}$} то очевидно, что расстояние между двумя его подмножествами может быть неограничено даже при фиксированной размерности пространства.
Если разрешить телу быть невыпуклым, то существует контрпример в виде <<ежа>>: тело, вытянутое вдоль каждой оси координат.
Поэтому видится целесообразным поставить следующую гипотезу

\begin{theorem}{Гипотеза 4.}
Пусть $\varepsilon$~--- данное число в интервале $(0,1)$, $K_n$~--- $n$-мерное выпуклое тело единичного объема инвариантное относительно произвольной перестановки координат.
Внутри $K_n$ выбраны два множества $A$ и $B$, каждое объема $\varepsilon$.
Тогда расстояние между $A$ и $B$ не больше, чем некоторая константа $D = D(\varepsilon)$, и не зависит от размерности пространства~$n$.
\end{theorem}

Например, можно рассмотреть случай правильных многомерных многогранников: тетраэдра и октаэдра.


\textbf{Благодарности.} Автор признателен своим научным руководителям А.\,Я.\,Белову и А.\,В.\,Михалеву за постановку задач и внимание к работе. Исследование поддержано грантом РФФИ № 14-01-00548.

\spisoklit
{\small\wrefdef{2}

\wref{1}
{Ширяев\,А.\,Н.}
Вероятность. МЦНМО, Москва, 2007

\wref{2}
{
Fomenko\,A.\,T.} Variational problems in topology, The geometry of length, area and volume, Gordon and Breach Science Publishers, New York, 1990, x+225 pp.

\wref{3}
{Дао Чонг Тхи, Фоменко\,А.\,Т.} Минимальные поверхности и проблема Плато, Наука, М., 1987, 312 с.  mathscinet  zmath; англ. пер.: Dao Trong Thi, A.\,T.\,Fomenko, Minimal Surfaces, Stratified Multivarifolds and the Plateau Problem, Translation of Mathematical Monographs, 84, Amer. Math. Soc., 1991, 404 с

\end{document}